\documentclass[a4paper]{article}

\begin{document}

\title{Generalized Projection Operators in\\Geometric Algebra}
\author{T.A. Bouma\\{\it University of Amsterdam}}
\date{April 15th, 2001 CE}
\maketitle

\begin{abstract}
Given an automorphism and an anti-automorphism of a semigroup of a
Geometric Algebra, then for each element of the semigroup a (generalized)
projection operator exists that is defined on the entire Geometric Algebra.  A
single fundamental theorem holds for all (generalized) projection operators.
This theorem makes previous projection operator formulas\cite{Hestenes-1}
equivalent to each other.  The class of generalized projection operators
includes the familiar subspace projection operation by implementing the
automorphism `grade involution' and the anti-automorphism `inverse' on the
semigroup of invertible versors.  This class of projection operators is studied
in some detail as the natural generalization of the subspace projection
operators.  Other generalized projection operators include projections onto
{\em any} invertible element, or a weighted projection onto {\em any} element.
This last projection operator even implies a possible projection operator for
the zero element.
\end{abstract}

\section{Introduction}
We introduce a class of generalized projection operators on a Geometric Algebra
$\mathcal{G}$ indexed by a nonempty set of generators, $G$, and two functions
from $G$ to $\mathcal{G}$ denoted: $A \mapsto \overline{A}$ and $A \mapsto A^{
\dagger}$ such that $G$ is closed under the geometric product, $\overline{AB}
= \overline{A}\hspace{0.35em}\overline{B}$, and $(AB)^{\dagger}=B^{\dagger}A^{
\dagger}$.  Each class of projection operators includes a function $P_A$ from
$\mathcal{G}$ to $\mathcal{G}$ for each element $A$ of $G$ defined by

\begin{equation}\label{definition}
P_A(X)= \frac{1}{2}(X-\overline{A}XA^{\dagger})
\end{equation}

The paper begins with the statement and proof of the Fundamental Theorem of
Projection Operators and then examines projection operators for specific choices
of $G$, $A \mapsto \overline{A}$, and $A \mapsto A^{\dagger}$.  The body of the
paper relates the fundamental theorem to familiar projection operators and
to novel projection operators.  The paper ends with a summary of future
work.

\section{Fundamental Theorem of Projection Operators}

The Fundamental Theorem of Projetion Operators (FToPO) states for any
projection operators $P_A$ and $P_B$

\begin{equation}\label{fundamental}
2P_A \circ P_B= P_A + P_B - P_{AB}
\end{equation}
Proof:
\begin{eqnarray*}
2P_A(P_B(X)) & = &\frac{1}{2}(X - \overline{B}XB^{\dagger} -
\overline{A}(X - \overline{B}XB^{\dagger})A^{\dagger})\\
2P_A(P_B(X)) & = & \frac{1}{2}(X - \overline{B}XB^{\dagger} -
\overline{A}XA^{\dagger} + \overline{A}(\overline{B}XB^{\dagger})A^{\dagger}) \\
2P_A(P_B(X)) & = & \frac{1}{2}(X - \overline{B}XB^{\dagger} -
\overline{A}XA^{\dagger} + \overline{(AB)}X(AB)^{\dagger}) \\
2P_A(P_B(X)) & = & \frac{1}{2}(X - \overline{B}XB^{\dagger} +
X -\overline{A}XA^{\dagger} -X + \overline{(AB)}X(AB)^{\dagger}) \\
2P_A(P_B(X)) & = &  P_B(X) + P_A(X) - P_{AB}(X)
\end{eqnarray*}

This theorem allows projection operators to be treated directly rather than
as derivative objects.  This can help for applications like those in
reference\cite{math.LA/0104102} where projection operators are used as 
fundamental objects of computation in Geometric Algebra.

\section{Familiar Projection Operators}
The most familiar projection is to let the set $G$ be the set of invertible
versors (the semigroup generated by the invertible blades), $\overline{A}$ be
the grade involution, and $A^{\dagger}$ be the inverse operation.
\subsection{Familiar Projections}

If $A$ is an invertible blade and $x$ is a vector, $P_A(x)$ is the projection
of $x$ onto the subspace characterized by $A$.

Proof:

\begin{eqnarray*}
P_A(x) & = & \frac{1}{2}(x-\overline{A}xA^{\dagger})\\
P_A(x) & = & \frac{1}{2}(xAA^{\dagger}-\overline{A}xA^{\dagger})\\
P_A(x) & = & \frac{1}{2}(xA-\overline{A}x)A^{\dagger}\\
P_A(x) & = & (x \rfloor A)A^{\dagger}
\end{eqnarray*}

So this class of projections is indeed familiar.

\subsection{Familiar Identities}

Here are seven formulas for the vector domain portion of the projection
operators $P_A$ and $P_B$ of two invertible blades $A$ and $B$ from
reference\cite{Hestenes-1}.

\begin{enumerate}
\item $AB = A \wedge B \Rightarrow P_B \circ P_A = 0$
\item $AB = A \wedge B \Rightarrow P_A \circ P_B = 0$
\item $AB = A \wedge B \Rightarrow P_{AB} = P_B + P_A$
\item $AB = A \rfloor B \Rightarrow P_{AB} = P_B - P_A$

\item $AB = A \rfloor B \Rightarrow P_A \circ P_B = P_A$
\item $AB = A \rfloor B \Rightarrow P_B \circ P_A = P_A$
\item $P_A \circ P_A = P_A$
\end{enumerate}

Using the FToPO it is easy to see that 1. $\Rightarrow$ 2. $\Rightarrow$ ...
$\Rightarrow$ 6. $\Rightarrow$ 7.

\begin{enumerate}
\item $\Rightarrow$ 2. from the FToPO and since $P_{AB} = P_{BA}$.  The latter
is clear since the reverse of an invertible blade is a nonzero scalar multiple
of itself.
\item $\Rightarrow$ 3. from the FToPO.
\item $\Rightarrow$ 4. because for invertible blades $A$ and $B$ $AB = A \rfloor
B \Rightarrow A^2B = A (A \rfloor B) = A \wedge (A \rfloor B)$ therefore
$P_B=P_{A^2B}=P_A + P_{A \rfloor B}$.
\item $\Rightarrow$ 5. from the FToPO.
\item $\Rightarrow$ 6. from the FToPO and since $P_{AB} = P_{BA}$.
\item $\Rightarrow$ 7. because for an invertible blade $A$, $AA=A \rfloor A$.
\end{enumerate}

\subsection{Versors and Blades}

The full generality of the FToPO interrelates the composition of projections
onto blades with projections onto versors.  Projecting onto a versor is a new
operation, but we will show a simple motivation, and that motivation will
reproduce the familiar projection onto a blade when the versor in question is,
in fact, a blade.  

Formula (\ref{definition}) clearly shows that $P_A(x)$ is the average of two
objects, namely $x$ and $-\overline{A}xA^{\dagger}$.  If $A$ is a versor then
the second object is always a vector.  Specifically, let $A=a_1a_2 ... a_r$ be
a versor.  Now define $x_0=x$ and inductively define $x_{i+1}= 
(-a_{r-i})x_i(a_{r-i}^{\dagger})$ then inductively it is clear that each
$x_i$ is a vector.  Expand $-\overline{A}xA^{\dagger}$ to get

\begin{eqnarray*}
-\overline{A}xA^{\dagger} & = & -(-a_1)(-a_2) ... (-a_r)x(a_r^{\dagger}) ...
(a_2^{\dagger})(a_1^{\dagger}) \\
&=&-(-a_1)(-a_2) ... (-a_rxa_r^{\dagger}) ... (a_2^{\dagger})(a_1^{\dagger})\\
&=&-(-a_1)(-a_2) ... (-a_rx_0a_r^{\dagger}) ... (a_2^{\dagger})(a_1^{\dagger})\\
&=&-(-a_1)(-a_2) ... (x_1) ... (a_2^{\dagger})(a_1^{\dagger})\\
&=&-(-a_1)(-a_2x_{r-2}a_2^{\dagger})(a_1^{\dagger})\\
&=&-(-a_1)(x_{r-1})(a_1^{\dagger})\\
&=&-(-a_1x_{r-1}a_1^{\dagger})\\
&=&-(x_r)
\end{eqnarray*}

As a family of versors $A(t)$
approaches the blade $B$, the vector $v(t)=\overline{A(t)}x(A(t))^{\dagger}$
becomes a vector whose rejection from $B$ remains the same as $x$, while the
projection of $v(t)$ swings around to become diametrically opposite the
projection of $x$.  Thus, $\frac{1}{2}(x-v(t))$ smoothly becomes the projection
of $x$ onto $B$.

\section{Novel Projection Operators}
The extension of projection to versors was required to fully utilize the FToPO
for blades and the interpretation of the projection of blades is truly an
explanation of the older, more familiar projection onto subspaces.  However,
there are also more formal extensions of the idea of projection, of which two
are explored here.

\subsection{Inverse Projection Operation}
The simplest formal extension of the familiar projection is to let the set $G$
be the set of all invertible elements, $\overline{A}$ be the grade involution,
and $A^{\dagger}$ be the inverse operation.  This clearly is just an
enlargement of the domain, $G$, of objects that can be projected onto.

We show that the interpretation of the projection onto a nonversor $W$ is
problematic.  This is because if $P_W(x)$ is a vector for each vector $x$, it
follows that $W$ is a versor.  Assume that for each vector $x$, $P_W(x)=
\frac{1}{2}(x-\overline{W}xW^{\dagger})$ is a vector.  Then isomorphically
embed the problem into a nontrivial, nondegenerate Geometric Algebra using a
LIFT as described in the appendix and define $f(x)=\overline{W}xW^{\dagger}$.
Clearly $f$ is a vector-valued linear function of a vector variable (i.e. $f$
is a linear transformation).  Furthermore $f$ is actually an orthogonal
transformation of the enlarged vector space.

\begin{eqnarray*}
(f(x))^2 & = & -f(x)\overline{f(x)}\\
& = & -(\overline{W}xW^{\dagger})\overline{(\overline{W}xW^{\dagger})}\\
& = & -(\overline{W}xW^{\dagger})(W\overline{x}\overline{W^{\dagger}})\\
& = & -\overline{W}x\overline{x}\overline{W^{\dagger}}=
\overline{W}x^2\overline{W^{\dagger}}\\
& = & x^2(\overline{W})(\overline{W^{\dagger}})=
x^2(\overline{WW^{\dagger}})=x^2\overline{1}\\
& = &x^2
\end{eqnarray*}

Since the enlarged space is nontrivial and nondegenerate
reference\cite{Hestenes-1} guarantees that there exists a nonzero versor $B$
that performs the same transformation, i.e. there exists a nonzero versor $B$
such that $\overline{B}xB^{\dagger}=\overline{W}xW^{\dagger}$ for each vector
$x$.  A short computation now shows that $x \rfloor (B^{\sim}W)=0$ for each
vector $x$, where $B^{\sim}$ denotes the reverse of $B$.  Note that $BB^{\sim}=
B^{\sim}B$ is a scalar and fix an arbitrary vector $x$.

\begin{eqnarray*}
\overline{B}xB^{\dagger} & = & \overline{W}xW^{\dagger}\\
\overline{B^{\sim}}(\overline{B}xB^{\dagger})(B^{\sim})^{\dagger} & = &
\overline{B^{\sim}}(\overline{W}xW^{\dagger})(B^{\sim})^{\dagger}\\
\overline{(B^{\sim}B)}x(B^{\sim}B)^{\dagger} & = &
\overline{(B^{\sim}W)}x(B^{\sim}W)^{\dagger}\\
x & = & \overline{(B^{\sim}W)}x(B^{\sim}W)^{\dagger}\\
x(B^{\sim}W) & = & \overline{(B^{\sim}W)}x\\
\frac{1}{2}(x(B^{\sim}W) - \overline{(B^{\sim}W)}x) & = & 0\\
x \rfloor (B^{\sim}W) & = & 0\\
\end{eqnarray*}

So by Lemma (\ref{lemma}) of the appendix, $\alpha=B^{\sim}W$ is a scalar.  Now 
$W=\frac{\alpha}{BB^{\sim}}B$ and since $W$ is a scalar multiple of $B$ it is a
versor too.

This means that projecting onto a nonversor, while defined, results in some
vectors going to nonvectors.  The interpretation of such a transformation is
an outstanding issue.

\subsection{Reverse Projection Operation}
The most general nontrivial projection operator is to let the set $G$ be the
set, $\mathcal{G}$ of all elements, $\overline{A}$ be the grade involution, and
$A^{\dagger}$ be the reverse.

As discussed in the previous section the interpretation of the projection onto
nonversors is problematic.  If $A$ is an invertible versor then
$\overline{A}xA^{\dagger}$ is proportional to $\overline{A}xA^{\-1}$, so in
the case where $A$ is an invertible versor the two projections are not very
different.  The inverse projection operation is the average, while the reverse
projection operation is the {\em weighted} average of $x$ and
$\overline{A}xA^{\dagger}$.

The class of reverse projection operators is defined for all multivectors, so
there is even a projection operator for the zero element, and in fact $P_O=
\frac{1}{2}$.  Projection operators for noninvertible elements are actually
quite interesting, for instance if an element, $A$, is idempotent ($A^2=A$) then
$P_A \circ P_A = \frac{1}{2}P_A$ follows from the FToPO easily.

The uses for the reverse projection operator are still unknown.  When
$A^{\dagger}= A^{-1}$ the reverse projection operator is the same as the
inverse projection operator.  When $A^{\dagger} \neq A^{-1}$ then the reverse
projection operator is a weighted average that depends on the scale of $A$.  So
possibly the reverse projection operator has use as a statistical projection
where the scale of $A$ determines the certainty of the element, or possibly the
actual {\em operation} on the elements will not be useful, but instead the
algebraic properties of the projection operator itself will provide a meaningful
(and useful) measure of the scale of a multivector.

\section{Conclusion}

This paper is a short introduction to a new class of projection operators in a
Geometric Algebra.  Even without taking projections onto new elements the
Fundamental Theorem of Projection Operators (FToPO) unifies and generalizes the
standard identities of projections onto subspaces.  The outright generalizations
of projection operators fall into three different potentially useful cases, each
of which calls for application or interpretation.  The first case is the
projection onto versors, which the author believes is the {\em natural}
generalization of the projection onto blades.  The second case is the projection
onto the zero element, which is simple enough that it can be appended
to any other class of projection operators and preserve the FToPO, and can
thereby introduce scale to projection operators.  The third case is the class of
weighted projection operators, $P_A$, that are sensitive to the scale of their
generators, $A$.

Each case is a call for further work.  The projection onto versors have a clear
interpretation and only need applications demonstrate its worth.  The zero
element was given a geometric interpretation in reference\cite{math.LA/0104102},
and that interpretation should be reconciled with the projection operator
presented here.  The weighted projection operators need both a solid
interpretation and applications and therefore will probably not be well
understood for some time to come.  Since the class of reverse projection
operators is a weighted projection operator that has a projection operator
$P_A$ for each element $A$ of the Geometric Algebra, there is at least the hope
that the reverse projection operators can help elucidate the geometric
properties of arbitrary elements of a Geometric Algebra.

\section*{Appendix}
The appendix contains two results used in the earlier proofs.

\subsection*{LIFT}
As taken from\cite{math.LA/0104102}, a LIFT (`linear injective function'
transformation) from one Geometric Algebra to another Geometric Algebra is
defined as a linear injective map that preserves the outer product and the
scalars.  In more detail, given two geometric algebras, $\mathcal{G}_1$ and
$\mathcal{G}_2$, and a linear injective function, $f$, from the vectors of
$\mathcal{G}_1$ to the vectors of $\mathcal{G}_2$ then $\underline{f}$ is a
LIFT between the two algebras, where $\underline{f}$ is the outermorphism of
$f$.

This paper uses a LIFT to isomorphically embed a degenerate Geometric Algebra
into a nondegenerate, nontrivial geometric algebra.  If the Geometic Algebra
is trivial it is just the scalars, and it is embedded into a Geometric Algebra
over a one-dimensional Euclidean vector space.  If the algebra is degenerate
then it is isomorphically embedded into a nondegenerate algebra as described
in reference\cite{math.LA/0104102}.

\subsection*{Contraction Lemma}
Let $\mathcal{G}$ be a Geometric Algebra over a nondegenerate finite dimensional
vector space $\mathcal{R}^{p,q}$ then

\begin{equation}\label{lemma}
x \rfloor A=0\hspace{1em}\forall x \in \mathcal{R}^{p,q} \Rightarrow A \in
\mathcal{R}
\end{equation}

Proof: Since $\mathcal{R}^{p,q}$ is nondegenerate it has an orthogonal basis of
invertible vectors $e_1, ... , e_n$.  The set $\{e_I: I \subset \{1, 2, ... , n
\}\}$ is a basis for $\mathcal{G}$ where $e_{\emptyset}=1$ and $e_I=\prod_{i_k
\in I} e_{i_k}$.  If $A = \sum \alpha^I e_I$ then $e_i \rfloor A = 0$ implies
that $\alpha^{I}=0$ when $i \in I$.  Since $e_i \rfloor A = 0$ for each $e_i$
it is clear that $\alpha^{I}=0$ for all $I \neq \emptyset$, therefore $A$ is a
scalar.

\end{document}